\documentclass{amsart}

\renewcommand{\phi}{\varphi}

\renewcommand{\span}{\mbox{span}\,}
\renewcommand{\r}{{\Bbb R}}
\newcommand{\z} {{\Bbb Z}}

\newcommand{\ddd}{,\dots,}
\newcommand{\lll}{\left(}
\newcommand{\rrr}{\right)}
\newcommand{\ex}[1]{e^{2\pi i{#1}}}
\newcommand{\exm}[1]{e^{-2\pi i{#1}}}

\newcommand{\cD}{{\mathcal D}}

\newcommand{\rk}{{\rm rk}}

\newcommand{\bQ}{{\mathbb Q}}

\newcommand{\bC}{{\mathbb C}}

\newcommand{\be}{\begin{equation}}
\newcommand{\ee}{\end{equation}}
\newcommand{\ba}{\begin{eqnarray}}
\newcommand{\ea}{\end{eqnarray}}
\newcommand{\ban}{\begin{eqnarray*}}
\newcommand{\ean}{\end{eqnarray*}}

\newtheorem{theorem}{Theorem}[section]

\theoremstyle{definition}
\newtheorem{definition}[theorem]{Definition}
\newtheorem{example}[theorem]{Example}

\theoremstyle{remark}

\theoremstyle{proposition}

\theoremstyle{proposition}

\newtheorem{theo}{Theorem}
\newtheorem{lem}[theo]{Lemma}
\newtheorem {prop} [theo] {Proposition}

\numberwithin{equation}{section}

\begin{document}

\title[On orthogonal $p$-adic  wavelet bases]
{On orthogonal $p$-adic  wavelet bases }

\thanks{ This research
was supported by Grant 12-01-00216 of RFBR and state-financed
project 9.38.62.2012 of SPbGU}

\author{S.~Evdokimov}
\address{St.Petersburg Department of Steklov Institute  of Mathematics, St.Petersburg,
 Fontanka-27, 191023 St.Petersburg, Russia  and
 St.Petersburg State University, Universitetskii pr.-35,
198504 St.Petersburg, Russia}
 \email{evdokim@pdmi.ras.ru}

\author{M.~Skopina}
\address{St.Petersburg State University, Universitetskii pr.-35,
198504 St.Petersburg, Russia}
\email{skopina@MS1167.spb.edu}
\thanks{}

.

\date{}


\begin{abstract}

A variety of different orthogonal wavelet bases has been found in $L_2(\r)$
for the last three decades. It appeared that similar constructions also exist for
functions defined on some other algebraic structures,
such as the Cantor and Vilenkin groups
and
local fields of positive characteristic.
In the present paper we show that the situation is quite  different
for the field of $p$-adic numbers. Namely, it is proved that any
orthogonal wavelet basis consisting of band-limited (periodic) functions is a
modification of Haar basis. This is a little bit unexpected
because from the wavelet theory point of view,
the additive group of $p$-adic numbers looks very similar to
the Vilenkin group where analogs of the Daubechies wavelets (and even band-limited ones) do exist.
We note that all  $p$-adic wavelet bases and frames appeared in the literature
consist of Schwartz-Bruhat functions (i.e., band-limited and compactly supported ones).
\end{abstract}

\maketitle

\section{Introduction}
\label{s1}

In the early nineties a general scheme for
constructing wavelets of real argument was developed. This
scheme is based on the concept of multiresolution analysis (MRA in the sequel)
introduced by Meyer and Mallat~\cite{Mallat-1}, \cite{Meyer-1} (see also,
~\cite{31}, ~\cite{NPS}). The theory allowed to construct orthogonal
wavelet bases  essentially different from the Haar basis, and, in particular, the Daubechies
wavelets \cite{31} which are actively implemented in signal processing and  other
engineering areas.

It appeared that similar constructions also exist for functions
defined on some other algebraic structures, such as the Cantor and
Vilenkin groups, local fields of positive and zero characteristic,
and adele rings (see, e.g.,~\cite{Lang}, \cite{PF}, \cite{F},
~\cite{Beh1}, \cite{Beh2}, \cite{S-Sk-1}, \cite{AES2}, \cite{E}).

 In the $p$-adic setting, the situation is as follows.
In 2002 the first $p$-adic wavelet basis for ${ L}_2(\bQ_p)$ where $\bQ_p$ is the $p$-adic field,
was found in \cite{Koz0}. It is an analog of the Haar basis for reals
and consists of translations by elements of the $p$-adic interval
$$
I_p=\Big\{\frac{k}{p^n}\in \bQ_p:\, k=0\ddd p^n-1,\ n=0,1\dots\Big\}
$$
and $p$-adic dilations, of $p-1$ functions (wavelet
functions in the sequel).
In~\cite{S-Sk-1} the notion of $p$-adic MRA was introduced
and a general scheme for its construction was described.
Following Meyer and Mallat one starts with a scaling function,
i.e. a function $\phi$ whose  $I_p$-translations
form an orthonormal system, and the sequence of spaces
$$
V_m=\overline{\span\Big\{\phi\lll\frac{x}{p^{m}}-a\rrr:\  a\in I_p\Big\}}, \quad m\in\z,
$$
is increasing and dense in ${ L}_2(\bQ_p)$.
The  scheme was used to construct the $p$-adic Haar MRA where
the characteristic function of the ring of $p$-adic integers $\z_p$ was taken as a scaling function.
This leads to the wavelet basis constructing in~\cite{Koz0}.

Next it was interesting to find other MRAs
in order to obtain essentially different wavelet bases.
However all these attempts were unsuccessful.  Although some other scaling functions
were found in~\cite{Kh-Sh-S}, it appeared that all of them
lead to the same Haar MRA. An explanation of this failure
was given in~\cite{AES2}, where the authors proved that the Haar MRA
is a unique MRA in the sense of the definition given in~\cite{S-Sk-1},
that is generated by a test (Schwartz-Bruhat)  scaling function.
So, it is not possible
to find essentially new orthogonal $p$-adic wavelet bases as MRA-based ones.
(It is interesting to note that as it was shown in \cite{ES1}, \cite{ES2}, there exist infinitely many
non-orthogonal MRAs in the sense of the definition given in~\cite{AES2}, i.e. those for which the $I_p$-translations
of the scaling function do not form an orthogonal system).

Now the following question arises: do there exist orthogonal wavelet bases
not generated by the Haar MRA?
In the present paper we answer to this question in the negative
restricting ourselves to the bases consisting of band-limited functions
(the class of band-limited functions contains the test functions, and
coincides with the class of periodic functions in ${ L}_2(\bQ_p)$).
Moreover, it is proved that any orthogonal wavelet basis consisting of test functions is
equivalent in a natural sense to the Haar basis constructed in \cite{Koz0}.

Let us discuss what the words "basis is generated by the Haar MRA" mean.
Following the standard scheme for constructing  MRA-based wavelets,
one defines the Haar wavelet spaces $W_m$  (the orthogonal complement
of $V_m$ in $V_{m+1}$) and find
a set of wavelet functions whose $I_p$-translations
form an orthonormal basis for $W_0$.
We call such a set standard;
the corresponding wavelet basis for $L_2(\bQ_p)$ consisting of $I_p$-translations
and $p$-adic dilations of its elements, is called a standard Haar basis.
Any standard set of wavelet functions consists of $p-1$ elements belonging to $W_0$.
Moreover, a set of functions that is unitary equivalent to
a standard set is standard too, but there are standard sets
which are not unitary equivalent to each other.
All standard sets consisting of test functions
were described in~\cite{Kh-Sh-S1}.

There exists, however, a lot of non-standard orthogonal
wavelet bases generated by the Haar MRA.
Indeed, it is easy to see
that each standard wavelet function $\psi$ can be replaced by $p$
functions (not belonging to $W_0$)
the set of $I_p$-translations and $p$-adic dilations of
which, coincides with that of $\psi$ (see Section 2 for details). So, the same standard
wavelet basis can be generated by a non-standard set of wavelet functions.
Moreover, if one applies a unitary transform to the latter set,
then a new orthogonal wavelet basis different from a standard one, can be obtained.
Of course these two operations (and their inverses)
can be repeated several times, and any basis obtained
in this way is actually generated by the Haar MRA. It is natural to call it
a ''damaged '' Haar basis.
We observe that all non-standard orthogonal $p$-adic wavelet bases we saw in the
literature (see, e.g.,~\cite{Ben-Ben}, ~\cite{Kh-Sh1}) are of that type.
In the present paper we  prove that any orthogonal $p$-adic wavelet basis whose
elements are band-limited is such a ''damaged '' Haar basis.

The paper is organized as follows. Notation and basic facts on $p$-adic analysis we need, are concentrated in Section~\ref{s2}.
Section~\ref{s3} contains basic definitions and statements of main results.
In Section~\ref{s4} we prove auxiliary results on general and periodic
vector-functions  generating orthonormal wavelet systems and especially orthonormal wavelet bases. Section~\ref{s5} contains the proofs
of Theorems~\ref{t1}, \ref{t3} and \ref{t2}.

\section{Notations and basic facts}
\label{s2}
Here and in what follows, we systematically use the
notation and results from~\cite{Vl-V-Z}.

As usual, by $\z$, $\bQ$ $\r$ and $\bC$ we denote the ring of rational integers and the fields of rational,
real and complex numbers, respectively.

The field $\bQ_p$ of $p$-adic numbers is the completion
of the field $\bQ$ with respect to the
$p$-adic norm $|\cdot|_p$ defined as follows:
$$
|x|_p=
\begin{cases}
p^{-\gamma},\ &\text{if }\ x=p^\gamma\frac m{n}\ne0;\\
0,\ &\text{if }\ x=0.
\end{cases}
$$
where $\gamma\in \z$
and $m$, $n$ are integers not divisible by $p$. The extension of this norm to $\bQ_p$ is also denoted by $|\cdot|_p$.

The norm $|\cdot|_p$ is non-Archimedean, i.e. satisfies the strong triangle inequality
$$
|x+y|_p\le \max(|x|_p,|y|_p),\quad x,y\in \bQ_p.
$$
Any $p$-adic number $x\ne 0$ can uniquely be written in the form
\begin{equation}
\label{2}
x=\sum_{j=\gamma}^\infty x_jp^j
\end{equation}
where $\gamma\in \z$ and $x_j\in\{0,1,\dots,p-1\}$ with $x_\gamma\ne 0$.
The fractional part $\{x\}_p$ of the number $x$
equals by definition $\sum_{j=\gamma}^{-1} x_jp^j$.
We also set $\{0\}_p=0$.

The ring of $p$-adic integers $\z_p$ and the $p$-adic interval $I_p$ are defined by
 $$
 \z_p=\{x\in \bQ_p:\, \{x\}_p=0\},\quad I_p=\{x\in \bQ_p:\, \{x\}_p=x\},
 $$
We observe that the translations of $\z_p$ by the elements of $I_p$ are mutually disjoint
and the union of them equals $\bQ_p$. Besides, $ \z_p=\{x\in \bQ_p:\, |x|_p\le1\}$.

The additive character $\chi_p$ of the field $\bQ_p$ is defined by
$$
\chi_p(x)=e^{2\pi i\{x\}_p},\quad x\in\bQ_p.
$$

The field $\bQ_p$ is locally compact. Denote by $dx$ the normalized Haar measure on it.
By definition this measure  is  positive, invariant under translations, i.e., $d(x+a)=dx$ for all $a\in\bQ_p$,
and satisfies the condition $\int_{\z_p}\,dx=1$. Moreover,
$$
d(ax)=|a|_p\, dx,\quad a\in\bQ_p\setminus\{0\}.
$$
The Hilbert space of all complex-valued functions on $\bQ_p$ summing with square with respect to
the measure  $dx$, is denoted by ${ L}_2(\bQ_p)$.

Denote by ${\cD}$ the $p$-adic Bruhat-Schwartz space, i.e. the linear space of locally-constant compactly
supported functions defined on $\bQ_p$ (so-called test functions).
This space is a $p$-adic analog of the  Schwartz space in the real analysis.

The Fourier transform of a function $f\in {\cD}$ is defined as
$$
{\widehat f}(\xi)=\int_{\bQ_p}\chi_p(\xi x)f(x)\,dx,
\ \ \ \xi \in \bQ_p,
$$
This yields a linear isomorphism taking ${\cD}$ onto
${\cD}$. It can uniquely be extended to a linear isomorphism of ${ L}_2(\bQ_p)$.
Moreover, the Plancherel equality holds
$$
\int\limits_{\bQ_p}f(x)\overline{g(x)}\,dx=
\int\limits_{\bQ_p}\widehat f(\xi)\overline{\widehat g(\xi)}\,d\xi,\quad f,g\in L^2(\bQ_p).
$$

A function $f\in { L}_2(\bQ_p)$ is said to be band-limited if its Fourier transform $\widehat f$
is compactly supported. Given $m\in\z$, the function  $\widehat f$ is supported on
$p^{-m}\z_p$ if and only if $f$ is $p^m$-periodic, i.e.
$$
f(x+p^m)=f(x),\quad x\in\bQ_p.
$$
Thus the spaces of band-limited and periodic functions coincide.
Moreover, the $p$-adic Bruhat-Schwartz space ${\cD}$ equals the space consisting
of all both band-limited (periodic) and compactly supported functions in  $L_2(\bQ_p)$.

Below for $f\in L_2(\bQ_p)$, $m\in\z$ and $a\in \bQ_p$, we set
\be
\label{0}
f_{m, a}(x)=p^{m/2}f\lll\frac{x}{p^{m}}-a\rrr,\quad x\in \bQ_p.
\ee

Let $\phi$ be a characteristic function of the set $\z_p$.
For every integer $m$, the functions $\phi_{m,a}$, $a\in I_p$,
form an orthonormal system.
The sampling spaces $V_m$ of the Haar MRA are defined by
\be
V_m=\overline{\span\{\phi_{m,a}:\ a\in I_p\}},\quad m\in\z.
\label{1}
\ee
The union of all  spaces $V_m$ is dense in ${ L}_2(\bQ_p)$, and
$V_m\subset V_{m+1}$ for all $m$  (see the p-adic MRA theory~\cite{S-Sk-1}).
Any function $\phi_{m,a}$ is $p^m$-periodic.
Moreover,
$$
V_m=\{f\in L_2(\bQ_p):\, f\ \text{is}\ p^m\text{-periodic}\}.
$$
The wavelet spaces $W_m$ are defined by
\be
 W_m= V_{m+1}\ominus V_m,\quad m\in \z,
\label{2}
\ee
and we have the following orthogonal decomposition
\begin{equation}
\label{3}
L_2(\bQ_p)={\bigoplus\limits_{m\in\z}W_m}.
\end{equation}

Finally, let us define a translation operator $T$ on $L_2(\bQ_p)$ by
\be
\label{T}
(Tf)(x)=f(x-1),\quad  x\in \bQ_p.
\ee
Then the operator $T$ is unitary and the spaces $V_j$ and $W_j$ are $T$-invariant.
Moreover, $Tf=f_{0,1}$ in accordance with (\ref{0}).

\section{Main results}
\label{s3}


Let $r$ be a positive integer and $\Psi=(\psi^{(1)}\ddd \psi^{(r)})^T$ where $\psi^{(\nu)}\in L_2(\bQ_p)$
for all $\nu=1\ddd r$.

\begin{definition}
\label{de3}
The system of functions
$$\{\psi^{(\nu)}_{m,a}:\  m\in\z,\, a\in I_p,\, \nu=1\ddd r\}$$
is called  {\em wavelet system} generated by the vector-function $\Psi$
(or by its components $\psi^{(\nu)}$ called {\em wavelets}).
If this  system is an orthonormal basis for $ L_2(\bQ_p)$,
one says that $\Psi$
generates an {\it orthonormal wavelet basis} (ONWB in the sequel).
\end{definition}

The number $r=\rk(\Psi)$ is called the {\em rank} of $\Psi$. The vector-function $\Psi$
is called  {\em eigen} if every component of $\Psi$ is an eigenfunction of the translation operator~(\ref{T}).

Two vector-functions
$\Psi$ and $\Psi^\prime$ are said to be {\it unitary equivalent} if there exists
a unitary matrix $U$ such that $\Psi=U\Psi^\prime$. Evidently, if $\Psi$
generates an ONWB and  $\Psi^\prime$ is  unitary equivalent to $\Psi$, then
$\Psi^\prime$ generates an ONWB.

\begin{definition}
\label{de1} \rm
A vector-function generating an ONWB is called {\em standard Haar} if it is of rank $p-1$
and all its components are in $W_0$. The ONWB is called  {\em standard Haar  basis} in this case.

\end{definition}

It is well known that the vector-function $\Theta=(\theta^{(1)},\ddd\theta^{(p-1)})^T$ where
$$
\theta^{(\nu)}(x)=\phi(x)\chi_p\lll\frac{\nu x}{p}\rrr, \quad \nu=1\ddd p-1,
$$
is standard Haar (see~\cite{S-Sk-1}). It is called the {\it basic Haar vector-function}.
Moreover, $\Theta$ is eigen because
\be
\label{basic}
\theta^{(\nu)}(x-1)=\chi_p\lll-\frac{\nu}{p}\rrr\theta^{(\nu)}(x).
\ee

Given a standard Haar vector-function, one can easily "damage" it
to obtain another vector-function (which is non-standard)
generating the same ONWB.
Using such a technique together with replacing some vector-functions by
unitary equivalent ones, new orthogonal wavelet bases can be obtained.
This is illustrated as follows.

\begin{example}
\label{ex1} \rm
Let $p=2$, $\psi=\theta^{(1)}$. Then $(\psi)$ can be treated as a standard Haar vector-function
consisting of one function $\psi$. Set
$$\psi^{(1)}(x)=\sqrt2\psi(x/2),\quad
\psi^{(2)}(x)=\sqrt2\psi((x-1)/2).
$$
Evidently,  the sets $\{\psi^{(\nu)}_{m,a}: m\in\z,\, a\in I_p,\, \nu=1,2\}$
and $\{\psi_{m,a}: m\in\z,\, a\in I_p\}$ coincide. So the wavelet system
generated by $\psi^{(1)}, \psi^{(2)}$ is a standard Haar
basis, but $(\psi^{(1)}, \psi^{(2)})^T$ is not a standard Haar vector-function.
Similarly, the functions  $\psi^{(1)}, \psi^{(2,1)}, \psi^{(2,2)}$,
where
$$
\psi^{(2,1)}(x)=\sqrt2\psi^{(2)}(x/2), \quad
\psi^{(2,2)}(x)=\sqrt2\psi^{(2)}((x-1)/2),
$$
generate the same standard Haar basis.
Set
$$
\widetilde\psi^{(1)}=\frac{1}{\sqrt 2}(\psi^{(1)}+\psi^{(2,1)}), \quad
\widetilde\psi^{(2)}=\frac{1}{\sqrt 2}(\psi^{(1)}-\psi^{(2,1)}), \quad
\widetilde\psi^{(3)}=\psi^{(2,2)}.
$$
These functions generate an ONWB which
is not a standard Haar basis, but the vector-function
$\widetilde\Psi=(\widetilde\psi^{(1)}, \widetilde\psi^{(2)}, \widetilde\psi^{(3)})^T$
is unitary equivalent to a vector-function generating
a standard Haar basis. Next set
$$
\widetilde\psi^{(1,1)}(x)=\sqrt2\widetilde\psi^{(1)}(x/2), \quad
\widetilde\psi^{(1,2)}(x)=\sqrt2\widetilde\psi^{(1)}((x-1)/2).
$$
It is not difficult to see that the vector-function
$\widetilde\Psi^\prime=
(\widetilde\psi^{(1,1)},\widetilde\psi^{(1,2)},
\widetilde\psi^{(2)}, \widetilde\psi^{(3)})^T$
generates an ONWB, but $\widetilde\Psi^\prime$
does not generate a standard Haar basis, and $\widetilde\Psi^\prime$ is not
unitary equivalent to a vector-function generating
a standard Haar basis.
\end{example}

\begin{definition}
\label{de2} \rm
Two vector-functions $\Psi$, $\Psi^\prime$  are said to be {\em wavelet equivalent } if
there exist  vector-functions $\Psi_{0}\ddd\Psi_{N}$ such that
$\Psi_{0}=\Psi$, $\Psi_{N}=\Psi^\prime$, and for every $j>0$
either $\Psi_{j}$ is unitary equivalent to $\Psi_{j-1}$
or $\Psi_{j}$ and $\Psi_{j-1}$ generate the same wavelet system (as a set).
\end{definition}

Evidently, the property of vector-functions to be wavelet equivalent is preserved
under changing the order of components in one of them.
So, we can use the term "wavelet equivalent" for two sets of functions,
or to say that a vector-function is wavelet equivalent to a set of
functions (which can be considered as the components of another vector-function).

We note also that the  wavelet equivalence relation  is reflexive,
symmetric and transitive. Moreover, if $\Psi$ generates an ONWB and  $\Psi^\prime$
is  unitary equivalent to $\Psi$, then $\Psi^\prime$ generates an ONWB.

It is easy to see that any standard Haar basis consists of periodic functions.
So, if a vector-function is wavelet equivalent to a standard Haar
vector-function, then all its components are periodic. It appears that
the converse statement is also true.

\begin{theo}
\label{t1}
Any periodic vector-function $\Psi$ generating an ONWB is wavelet
equivalent to a standard Haar vector-function which can be taken eigen. If, moreover, $\Psi$ is compactly  supported, then
it is wavelet equivalent to the basic Haar vector-function ~$\Theta$.
\end{theo}

As a byproduct of Theorem~\ref{t1} we obtain the following statement.

\begin{theo}
\label{t3}
The rank $r$ of a periodic vector-function
generating an ONWB is divisible by $p-1$, in particular, $r\ge p-1$.
\end{theo}

When considering Example~\ref{ex1}, we see how to find ONWBs which are
not standard Haar bases. Every vector-function constructed
in this way is wavelet equivalent to a standard Haar one.
However, really the situation is more specific.
Namely,
at each step, where the vector-functions $\Psi_{j}$ and $\Psi_{j-1}$
generate the same wavelet basis, in fact $\Psi_{j}$ is obtained from
$\Psi_{j-1}$ by changing one of its components for two functions.
To describe all bases that can be constructed in such a way for an arbitrary $p$,
we say that $\Psi^\prime$ {\em is reducible} to $\Psi$ if
there exist  vector-functions $\Psi_{0}\ddd\Psi_{N}$ such that
$\Psi_{0}=\Psi$, $\Psi_{N}=\Psi^\prime$, and for every $j>0$,
either $\Psi_{j}$ is unitary equivalent to $\Psi_{j-1}$
or $\Psi_{j-1}$ is obtained from $\Psi_{j}$ by changing one of
its components, say $\psi_j^{(\nu)}(x)$, for $p$ functions
$\psi_{j}^{(\nu)}\big(\frac{x-k}p\big)$, $k=0\ddd p-1$.

It would be attractive to replace the words "wavelet equivalent"
by "reducible" in Theorem~\ref{t1}. However this is impossible
due to the following statement.

\begin{theo}
\label{t2}
There exists a vector-function  generating an ONWB that
cannot be reduced to a standard Haar vector-function.
\end{theo}

\section{Auxiliary results}
\label{s4}

\begin{lem}
\label{l1}
Let $m\ge0$ be an integer and $\psi\in L_2(\bQ_p)$. If $\psi=\sum_{a\in I_p}c_a\phi_{m, a}$ with $c_a\in\bC$,
then
\be
\sum_{\nu=0}^{p^m-1}\sum_{j=0}^{\infty}\sum_{b\in I_p}
|\langle f^{(\nu)}, \psi_{-j,b}\rangle|^2=\frac p{p-1}\sum_{a\in I_p}|c_a|^2
\label{6}
\ee
where $f^{(\nu)}=\phi_{m,\nu/p^{m}}$.
\end{lem}

{\bf Proof}. Set
$$A_{0}=\{0\}\cup\Big(I_p\cap\big(1-\frac{1}{p^m}, 1\big)\Big),\quad
A_\nu=I_p\cap\big(\frac{\nu-1}{p^m}, \frac{\nu}{p^m}\big],\ \
\nu=1\ddd p^{m}-1.$$
Then evidently, $I_p=\bigcup_{\nu=0}^{p^m-1}A_\nu$ and
$A_\nu\cap A_\mu=\emptyset$ for $\nu\ne\mu$.

Let $\nu\ne0$, $j\ge0$ and $b\in I_p$.
Then
\ban
\langle f^{(\nu)}, \psi_{-j,b}\rangle=
\sum\limits_{a\in I_p}c_a\int\limits_{\bQ_p}
p^{m/2}\phi\lll\frac x{p^m}-\frac{\nu}{p^m}\rrr
p^{(m-j)/2}\phi\lll\frac x{p^{m-j}}-\frac b{p^m}-a\rrr\,dx=
\\
p^{-j/2}\sum\limits_{a\in I_p}c_a\int\limits_{\bQ_p}\phi(x)
\phi\lll p^j x-\lll a+\frac b{p^m}-\frac \nu{p^m}\rrr\rrr\,dx=
p^{-j/2}c_{(\nu-b)/{p^m}}.
\ean
Since ${(\nu-b)/{p^m}}\in A_\nu
$ for all $b\in I_p$, and
each $a\in A_\nu$ can be represented in this form,
we have
\be
\sum_{b\in I_p}
|\langle f^{(\nu)}, \psi_{-j,b}\rangle|^2=p^{-j}\sum_{a\in A_\nu}
|c_a|^2.
\label{4}
\ee
Now let $\nu=0$. Then similarly, for $j\ge0$ we have
$\langle f^{(0)}, \psi_{-j,0}\rangle=p^{-j/2}c_{0}$ and
$\langle f^{(0)}, \psi_{-j,b}\rangle=p^{-j/2}c_{1-\frac b{p^m}}$, if $b\in I_p\setminus\{0\}$.
This yields
\be
\sum_{b\in I_p}
|\langle f^{(0)}, \psi_{-j,b}\rangle|^2=p^{-j}\sum_{a\in A_0}
|c_a|^2.
\label{5}
\ee
Adding~({\ref4}) and~({\ref5}) for all  $\nu=0\ddd p^m-1$ and $j=0,1,\dots$
we obtain~({\ref6}).$\Diamond$

\begin{prop}
\label{p1}
Let   functions $\psi^{(1)}\ddd \psi^{(r)}\in V_m$, $m\ge0$, generate
an orthonormal wavelet system. Then
\begin{itemize}
\item [(1)]
  $r\le(p-1)p^{m-1}$; in particular,
the wavelet system generated by $\psi\in V_0$ cannot be orthogonal.
\item [(2)]
If all functions $\psi^{(1)}\ddd \psi^{(r)}$ are in $W_{m-1}$,  $m>0$,
and generate an ONWB, then $r=(p-1)p^{m-1}$.
\end{itemize}
\end{prop}

{\bf Proof}.
Any function
$\psi^{(\nu)}$
is in $V_m$, so it can be
expanded as $\psi^{(\nu)}=\sum\limits_{a\in I_p}c^{(\nu)}_a\phi_{m a}$.
Let  $f^{(\mu)}$
be the function from Lemma~\ref{l1}.
Then using the Bessel inequality and~(\ref{6}), we have
\ba
p^m=\sum\limits_{\mu=0}^{p^m-1}\|f^{(\mu)}\|^2\ge
\sum_{\mu=0}^{p^m-1}\sum_{\nu=0}^{r}\sum_{j\in\z}\sum_{b\in I_p}
|\langle f^{(\mu)}, \psi^{(\nu)}_{-j,b}\rangle|^2\ge\hspace{2cm}
\label{7}
\\
\sum_{\nu=0}^{r}\sum_{\mu=0}^{p^m-1}\sum_{j=0}^\infty\sum_{b\in I_p}
|\langle f^{(\mu)}, \psi^{(\nu)}_{-j,b}\rangle|^2=
\sum_{\nu=0}^{r}\frac p{p-1}\sum_{a\in I_p}|c^{(\nu)}_a|^2=
\sum_{\nu=0}^{r}\frac p{p-1}\|\psi^{(\nu)}\|^2=\frac {rp}{p-1},
\nonumber
\ea
which
proves statement (1).

 Now if $\psi^{(\nu)}\in W_{m-1}$, then $\psi^{(\nu)}_{j,b}$ is orthogonal to $V_m$ for
all $j>0$, $b\in I_p$. Hence,
$\langle f^{(\mu)}, \psi^{(\nu)}_{-j,b}\rangle=0$ whenever $j<0$. It follows that
the second inequality in~(\ref7) can be replaced
by  equality. If, moreover, the functions $\psi^{(1)}\ddd \psi^{(r)}$
generate an ONWB, then due to the Parseval equality the first inequality in~(\ref7)
also can be replaced by  equality. This proves statement (2). $\Diamond$

To formulate Proposition~\ref{p7} which is a driver of the Theorem~\ref{t1} proof,
we need the following simple observation to be also used several times in what follows.

\begin{lem}
\label{l2}
Let $m, n$ be positive integers, $m\ge n$, and $f\in V_{m-n}$.
Then there exist functions
$f^{(k)}\in V_{m}$, $k=0\ddd p^n-1$,
the wavelet system generated by which coincides with
that generated by $f$. The functions $f^{(k)}$ can be given by
\be
f^{(k)}(x):=f_{n,k/p^n}(x)=p^{n/2} f\lll\frac {x-k}{p^n}\rrr,\quad k=0\ddd p^n-1.
\label{10}
\ee
\end{lem}

{\bf Proof.} The functions defined in~(\ref{10}) are what we need
because
$$
f^{(k)}\lll\frac x{p^j}-b\rrr=p^{n/2}f\lll\frac x{p^{j+n}}-\frac {k+b}{p^n}\rrr
$$
for all $j\in\z$ and $b\in\bQ_p$, and
$$I_p= \bigcup\limits_{k=0}^{p^{n}-1}
\left\{\frac {k+b}{p^n}:\, b\in I_p\right\}.\Diamond
$$



Below by a {\it non-trivial} linear combination of functions $\psi^{(1)}\ddd\psi^{(r)}$ we mean
any function $\alpha_1\psi^{(1)}+\dots + \alpha_r\psi^{(r)}$ such that $\alpha_{\nu}\ne0$ for some $\nu=1\ddd r$.

\begin{prop}
\label{p7} Let a  vector-function $\Psi$ of rank $r$ with all components in $V_m$ where  $m>0$,
generate an ONWB. Suppose that there exists a non-trivial linear
combination of the components
which is in $V_{m-1}$.
Then $\Psi$ is wavelet equivalent to a vector-function
of rank $r+p-1$, all components of which are in $V_m$.
 \end{prop}

 {\bf Proof.}
Let $\Psi=(\psi^{(1)}\ddd\psi^{(r)})^T$ and let $\alpha_1\psi^{(1)}+\dots + \alpha_r\psi^{(r)}$
be a non-trivial linear combination belonging to
$V_{m-1}$. Set $\Psi'=U\Psi$, where  $U$ is a unitary matrix whose
first row is the normalized vector $(\alpha_1\ddd \alpha_r)$. The
vector-function $\Psi'$ is unitary equivalent to $\Psi$, which
yields that $\rk(\Psi')=r$ and $\Psi'$ generates an ONWB. Since the first component
of $\Psi'$ is in $V_{m-1}$ and the other components are in $V_m$,
by Lemma~\ref{l2} there exists a vector function $\widetilde\Psi$
with all components in $V_m$ such that $\rk(\widetilde\Psi)=r+p-1$ and the wavelet systems generated by
$\Psi'$ and $\widetilde\Psi$ coincide.  $\Diamond$

\begin{prop}
\label{p8}
If a periodic vector-function $\Psi=(\psi^{(1)}\ddd \psi^{(r)})^T$ generates an ONWB,
then
\be
\Psi(x-1)=\sum\limits_{j=-\infty}^nA_j\Psi(p^{-j}x),\quad  x\in\bQ_p,
\label{12}
\ee
where $n\ge0$ and $A_j$ is an $r\times r$ matrix with complex entries.
\end{prop}

 {\bf Proof.}  The functions $\psi^{(\nu)}_{ja}$, $a\in I_p$, $j\in\z$, $\nu=1\ddd r$,
 form an orthonormal basis for $L_2(\bQ_p)$. Therefore,
for every $\mu$ we have
 $$
 \psi^{(\mu)}(x-1)=\psi^{(\mu)}_{01}(x)=\sum_{j\in\z}\sum_{a\in I_p}\sum_{\nu=1}^r
 \langle \psi^{(\mu)}_{01}, \psi^{(\nu)}_{ja} \rangle
 \psi^{(\nu)}_{ja}(x).
 $$
If $j\in\z$ and $a\in I_p\setminus\{0\}$,  then  there exists a
rational number $b=b(j,a)$,
such that $1-b\in I_p$, $a-p^{-j}b\in I_p$, whence
\ban
\langle \psi^{(\mu)}_{01}, \psi^{(\nu)}_{ja} \rangle=p^{j/2}
\int\limits_{\bQ_p} \psi^{(\mu)}(x-1)\overline{\psi^{(\nu)}(p^{-j}x-a)}\ dx=
\\
p^{j/2}\int\limits_{\bQ_p} \psi^{(\mu)}(x-(1-b))
\overline{\psi^{(\nu)}(p^{-j}x-(a-p^{-j}b))}\ dx=0,
\ean
i.e.,
\be
 \psi^{(\mu)}(x-1)=\sum_{j\in\z}\sum_{\nu=1}^r
 \langle \psi^{(\mu)}_{0,1}, \psi^{(\nu)}_{j,0} \rangle
 \psi^{(\nu)}_{j,0}(x).
 \label{11}
 \ee
 Since $\Psi$ is a periodic vector-function, there exists $m\in\z$ such that
 $\psi^{(\nu)}\in V_m$ for all $\nu$. Due to
 Proposition~\ref{p1}, the function  $\psi^{(\nu)}$ is not in $V_0$, i.e.
 $m>0$.  If $j\ge m$, then
 $p^{m-j}-p^{-j}\in I_p\setminus\{0\}$, which yields
\ban
\langle \psi^{(\mu)}_{0,1}, \psi^{(\nu)}_{j,0} \rangle=
p^{j/2}\int\limits_{\bQ_p} \psi^{(\mu)}(x+p^m-1)\overline{\psi^{(\nu)}(p^{-j}x)}\ dx=
\\
p^{j/2}\int\limits_{\bQ_p} \psi^{(\mu)}(x)
\overline{\psi^{(\nu)}(p^{-j}x-(p^{m-j}-p^{-j}))}\ dx=0.
\ean
 Substituting this into~(\ref{11}), we get~(\ref{12})  with $n=m-1$.
   $\Diamond$

\medskip 
Let $m>0$ be an integer. For any integer $l$ set
\be
V_{m,l}=\{f\in V_m:\, Tf=\exm{\frac{l}{p^{m}}}f\,\}
\label{18}
\ee
where $T$ is the translation operator~(\ref{T}). Thus, $V_{m,l}$ is the eigenspace of $T$ that corresponds to the
eigenvalue $\lambda_l=\exm{\frac{l}{p^{m}}}$. Obviously, $V_{m.0}=V_0$ and $V_{m,l}=V_{m,l'}$
whenever $l\equiv l'(\hspace{-2mm}\mod p^m)$. Moreover, since eigenfunctions belonging to different
eigenvalues are orthogonal, the spaces $V_{m,l}$ where  $0\le l\le p^m-1$, are pairwise orthogonal.

\begin{lem}
\label{p11}
Let $m>0$ be an integer. Then

\begin{itemize}

\item[(1)]
$V_{m,l}=V_{m+k,lp^k}$ for all $l$ and all $k\ge0$;

\item[(2)]
$W_{m-1}=\bigoplus\nolimits_{l\in S_m}V_{m,l}$ where $S_m=\{l\in\{0\ddd  p^m-1\}:\, p\ \text{does not divide}\ l\}$;

\item[(3)]
If $l\in S_m$, then the space $V_{m,l}$ is orthogonal to the wavelet system generated by any function
belonging to $V_{m,l'}$ with $l'\in S_m\setminus\{l\}$.
\end{itemize}
\end{lem}

{\bf Proof.} Statement (1) is obvious. Further, since the operator $T$
is unitary and $T^{p^m}=I$ on $V_m$, the space $V_m$ can be decomposed into the orthogonal sum
of its eigenspaces (see e.g. \cite[Ch.1]{S}):
$$
V_m=\bigoplus_{l=0}^{p^m-1}V_{m,l}.
$$
(Indeed,
if $f\in V_m$, then
$f=\sum_{l
}f_l$ where $f_l:=p^{-m}\sum_{j
}\ex{\frac{j}{p^{m}}}T^jf$ is in $V_{m,l}$).
Moreover, by statement (1) any space $V_{m.l}$ with $l$ coprime to $p$ is orthogonal to the space
$V_{m-1}=\bigoplus_{l=0}^{p^{m-1}-1}V_{m,lp}$
whereas any space $V_{m.l}$ with $l$ dividing $p$ is a subspace of $V_{m-1}$. Thus statement (2)
follows from the definition of $W_{m-1}$. To prove statement (3) it suffices to
note that by (\ref{3}) the spaces $W_j$, $j\in\z$, are pairwise orthogonal and use statement (2).
$\Diamond$

\begin{prop}
\label{p10}
Let a  vector-function $\Psi$
generate
an ONWB. Suppose that for some $m>0$ any non-trivial linear combination
of its components is in $V_m\setminus V_{m-1}$.
Then $r(\Psi)=(p-1)p^m$ and $\Psi$ is unitary equivalent to an eigen vector-function
every component of which belongs to $W_{m-1}$.

\end{prop}

 {\bf Proof.}
 Due to Proposition~\ref{p8},
identity~(\ref{12}) holds with some $n\ge0$. Assume that $n>0$ and
$A_{n}$ is a non-zero matrix. Then by the hypothesis at least one
component of the vector-function $A_{n}\Psi(p^{-n}x)$ is in
$V_{m+n}\setminus V_{m+n-1}$, whereas all components of the left-hand
side of~(\ref{12}) and of the sum $\sum_{j=-\infty}^{n-1}A_j\Psi(p^{-j}x)$ are in $V_{m+n-1}$,
a contradiction. Therefore, $n=0$, i.e.,
\be \Psi(x-1)=\sum_{j=-\infty}^0 A_j\Psi(p^{-j}x),\quad x\in\bQ_p, \label{14}
\ee
Iterating this identity $p^m$ times and
taking into account that $\Psi(x-p^m)=\Psi(x)$ for all
$x\in\bQ_p$, we obtain that every component of the vector-function
$(I-A_0^{p^m})\Psi$ belongs to $V_{m-1}$,
where $I$ is the $r\times r$ identity matrix with $r=\rk(\Psi)$.
But if $I\ne A_0^{p^m}$, then by the hypothesis at least one
component of this vector-function
is in $V_m\setminus V_{m-1}$, a contradiction.
Thus, \be A_0^{p^m}=I.\label{13'}\ee

It follows from~(\ref{14}) and the Parseval equality that the
Euclidean norm of each row of $A_0$ does not exceed  $1$. On the
other hand, by the Adamar inequality the product of these norms is
not less than $|\det A_0|$ which equals 1 due to (\ref{13'}).
Thus, $A_0$ is a unitary matrix. Again using the Parseval
equality, we obtain that $A_j=0$ for all $j<0$, i.e. \be
\Psi(x-1)=A_0\Psi(x),\quad  x\in\bQ_p. \label{15} \ee

Let $\lambda_1\ddd\lambda_r$ be eigenvalues of $A_0$, and let $D$
be the $r\times r$ diagonal matrix with $\lambda_1\ddd\lambda_r$ on the diagonal.
There exists a unitary matrix $U$ such that $A_0=UDU^{-1}$. Set
$\widetilde\Psi=U\Psi$, and rewrite~(\ref{15}) as
$$
\widetilde\Psi(x-1)=D\widetilde\Psi(x), \quad  x\in\bQ_p.
$$
So every component of $\widetilde\Psi$ is an eigenfunction of $T$. It cannot belong to $V_{m-1}$
by the hypothesis, thus it is in $W_{m-1}$ by statement (2) of Lemma~\ref{p11}.
$\Diamond$

\begin{prop}
\label{p11'}
Let $\Psi$ be an eigen vector-function generating an ONWB. Suppose that $r(\Psi)=(p-1)p^{m-1}$ and every component of $\Psi$
belongs to $W_{m-1}$.
Then $\Psi$ is wavelet equivalent to an eigen standard Haar vector-function.
\end{prop}
{\bf Proof.} Let $\Psi=(\psi^{(1)}\ddd \psi^{(M)})^T$ where $M=(p-1)p^{m-1}$. By the hypothesis and
statement (2) of Lemma~\ref{p11} we have $\psi^{(\nu)}\in V_{m,l_\nu}$ for all $\nu=1\ddd M$
where ${l_\nu}\in S_m$.
If the mapping
$$
F:\{1\ddd M\}\to S_m,\quad \nu\mapsto l_\nu
$$
is not a surjection, then there exists $l\in S_m$ not belonging to the image of $F$. So by
statement~(3) of Lemma~\ref{p11}, any function in the space $V_{m,l}\ne\{0\}$ is orthogonal to the wavelet system generated by $\Psi$
which is a basis for $L_2(\bQ_p)$ by the hypothesis, a contradiction. Since $\#S_m=M$ we conclude that $F$ is a bijection.

For $\mu=1\ddd p-1$ set
$$
f^{(\mu)}(x)=p^{1-m}\sum_{\nu\in T_\mu}\psi^{(\nu)}(p^{m-1}x),\quad x\in \bQ_p,
$$
where
$$
T_\mu=F^{-1}(S_{m,\mu})\quad \text{with}\quad S_{m,\mu}=\{l\in S_m:\,l\equiv \mu(\hspace{-3mm}\mod p)\}.
$$
We have $\|f^{(\mu)}\|=1$
for all $\mu$ because $\# T_\mu=\# S_{m,\mu}=p^{m-1}$. Moreover,
\ba
\nonumber f^{(\mu)}(x-1)=c\sum_{\nu\in T_\mu}\psi^{(\nu)}(p^{m-1}x-p^{m-1})=c\sum_{\nu\in T_\mu}\exm{\frac{l_\nu}{p}}\psi^{(\nu)}(p^{m-1}x)=\\
c\sum_{\nu\in T_\mu}\exm{\frac{\mu}{p}}\psi^{(\nu)}(p^{m-1}x)=\exm{\frac{\mu}{p}}f^{(\mu)}(x)
\label{19}
\ea
where $c=p^{1-m}$.

By Lemma~\ref{l2} the wavelet system generated by a function $f^{(\mu)}$ coincides with that generated by the functions
$f_{m-1,k/p^{m-1}}^{(\mu)}$, $k=0\ddd p^{m-1}-1$, where
$$
f_{m-1,k/p^{m-1}}^{(\mu)}(x)=p^{(m-1)/2}f^{(\mu)}\lll\frac{x-k}{p^{m-1}}\rrr.
$$
On the other hand, since $S_{m,\mu}=\{l=\mu+pj:\, j=0\ddd p^{m-1}-1\}$ we have
\ban
f_{m-1,k/p^{m-1}}^{(\mu)}(x)=c^{1/2}\sum_{\nu\in T_\mu}\psi^{(\nu)}(x-k)=c^{1/2}\sum_{\nu\in T_\mu}\exm{\frac{l_\nu k}{p^m}}\psi^{(\nu)}(x)=\\
c^{1/2}\sum_{l\in S_{m,\mu}}\exm{\frac{lk}{p^m}}\psi^{(\nu_l)}(x)=c^{1/2}\sum_{j=0}^{p^{m-1}-1}\exm{\frac{(\mu+pj)k}{p^m}}\psi^{(\nu_{j,\mu})}(x)=\\
\exm{\frac{\mu k}{p^m}}\sum_{j=0}^{p^{m-1}-1}c^{1/2}\exm{\frac{jk}{p^{m-1}}}\psi^{(\nu_{j,\mu})}(x)
\ean
where $\nu_{j,\mu}=\nu_l=F^{-1}(l)$ for $l=\mu+pj$. However it is well known that the matrix
$$\lll p^{-n/2}\exm{\frac{jk}{p^{n}}}\rrr_{j,k=0}^{p^{n}-1},$$
is unitary for any integer $n\ge0$.
So the system of functions
$$\{f_{m-1,k/p^{m-1}}^{(\mu)}: k=0\ddd p^{m-1}-1\}$$ is unitary equivalent to the system
$\{\psi^{(\nu_{j,\mu})}:\, j=0\ddd p^{m-1}-1\}=\{\psi^{(\nu)}: \nu\in T_\mu\}$.
Since the sets $ T_\mu$ are pairwise disjoint and the union of them equals $\{1\ddd M\}$, we conclude that the vector-function
$\Psi'=(f^{(1)}\ddd f^{(p-1)})$ is wavelet equivalent to~$\Psi$. The components of $\Psi$ are in $W_{m-1}$,
so the components of $\Psi'$ are in $W_0$. Thus, $\Psi'$ is a standard Haar vector function. It is eigen due to (\ref{19}).

\section{Proof of main results}
\label{s5}

{\bf Proof of Theorems~\ref{t1} and \ref{t3}.}
In what follows, we say that a vector-function is in $V_j$ if all its components
are in $V_j$.
Let $\Psi=\Psi_{0}$ be a $p^m$-periodic vector-function generating an ONWB.
Then $\Psi_{0}$ is in $V_m$ and $m>0$ because of Proposition~\ref{p1}.
Assume that some non-trivial linear
combination of components of $\Psi_{0}$ is in $V_{m-1}$.
Due to Proposition~\ref{p7},  $\Psi_{0}$ is wavelet equivalent to a vector-function
$\Psi_{1}$ in $V_m$ such that $\rk(\Psi_{1})=\rk(\Psi_{0})+p-1$.
Similarly, if some non-trivial linear
combination of components of $\Psi_{1}$  is in $V_{m-1}$, then
 $\Psi_{1}$ is wavelet equivalent to a vector-function
$\Psi_{2}$ in $V_m$ such that $\rk(\Psi_{2})=\rk(\Psi_{1})+p-1$.
Continue this process while possible. At each step we obtain a new vector-function
belonging to $ V_m$ and generating an ONWB.
The rank of the vector-functions strictly increases.
However by Proposition~\ref{p1}, it
does not exceed  $(p-1)p^{m-1}$. Hence the
process will stop after a finite number of steps, say $N$, i.e.,
$\Psi$ is wavelet equivalent to a vector-function
$\Psi_{N}$  such that any non-trivial linear combination
of its components is in $V_m\setminus V_{m-1}$
(in particular, it may happen that $N=0$). By Proposition~\ref{p10}, $\rk(\Psi_{N})=(p-1)p^{m-1}$. So
$$
\rk(\Psi)=\rk(\Psi_{0})\equiv\rk(\Psi_{1})\equiv\dots\equiv\rk(\Psi_{N})\equiv0(\hspace{-3mm}\mod p-1)
$$
whence Theorem~\ref{t3} follows. Besides, Propositions~\ref{p10} and \ref{p11'} imply that
$\Psi_{N}$ (and hence $\Psi)$ is wavelet equivalent to an eigen standard Haar vector-function,
which proves the first part of Theorem~\ref{t1}.

Now suppose that every component of $\Psi$ is compactly supported.
By above without loss of generality we can also assume that $\Psi=(\psi^{(1)}\ddd \psi^{(p-1)})^T$
is a standard Haar vector  function.
Every function $\psi^{(\mu)}$ is in $W_0$. So it can be expanded on the basis $\theta^{(\nu)}_{0, a}$, $a\in I_p$,
 $\nu=1\ddd p-1$, of this space where $\theta^{(\nu)}$ is the $\nu$-th component of the basic Haar vector-function $\Theta$.
Since $\psi^{(\mu)}$ is compactly supported, the expansion is finite. Therefore, there exists an integer
$n\ge0$ such that
\begin{equation}
\label{101**}
\psi^{(\mu)}=\sum_{\nu=1}^{p-1}\sum_{k=0}^{p^n-1}c_{\nu, k}^{\mu}
\theta_{0, k/p^n}^{(\nu)},
\quad \mu=1,2,\dots,p-1,
\end{equation}
where
$c_{\nu, k}^\mu$ is a complex number.

Let us denote by $\Psi'$ (resp. $\Theta'$) the vector-function of rank $(p-1)p^n$
with components numerated by pairs $(\nu,k)$ where $\nu=1\ddd p-1$ and $k=0\ddd p^n-1$,
such that the $(\nu,k)$-th component
equals $\theta^{(\nu)}_{n, k/p^n}$ (resp. $\psi^{(\nu)}_{n, k/p^n}$).
Due to Lemma~\ref{l2}, the vector-functions $\Psi'$ and
$\Psi$ (as well as $\Theta'$ and $\Theta$) generate the same ONWB.
So, it only remains to check that $\Psi'$ and $\Theta'$ are unitary equivalent.

To do this, from~(\ref{101**}) we derive that
$$
\psi_{n, l/p^n}^{(\mu)}=\sum_{\nu=1}^{p-1}\sum_{k=0}^{p^n-1}c_{\nu k}^{\mu}
\theta_{n,(k+l)/p^n}^{(\nu)},\quad l=0\ddd p^n-1,
$$
and taking into account that $\theta_{n,(a+p^n)/p^n}^{(\nu)}=
\chi_p\big(-\frac{\nu}{p}\big)\theta_{n,i/p^n}^{(\nu)}$ for all $a\in\bQ_p$ due to (\ref{basic}),
we conclude that each function $\psi^{(\mu)}_{n, l/p^n}$
is a linear combination of the functions $\theta^{(\nu)}_{n, k/p^n}$,
$\nu=1\ddd p-1$, $k=0\ddd p^n-1$. Hence,
$$
\Psi'=U\Theta'
$$
where $U$ is a matrix with complex entries. But the matrix $U$ is unitary because the components of each of the vector-functions
$\Psi'$ and $\Theta'$ form an orthonormal system of the same rank.
 $\Diamond$

{\bf Proof of Theorem~\ref{t2}.}
Let $p=2$, and $\theta=\theta^{(1)}$.
For $k=0,1$ set
$$
f^{(k)}(x)=\frac{\sqrt2}{2}\sum_{l=0}^1\exm{\frac{l(1+2k)}{4}}\theta^{}_{1,l/2}(x)=\sum_{l=0}^1\exm{\frac{l(1+2k)}{4}}
\theta^{}\lll \frac {x-l}{2}\rrr.
$$
The vector-function ($f^{(0)}, f^{(1)}$) generates an ONWB because it is unitary equivalent
to ($\theta_{1,0}$, $\theta_{1,1/2}$) and by Lemma~\ref{l2} the latter generates the same ONWB as the basic
vector-function $\Theta=(\theta)$.
Note that $f^{(0)}, f^{(1)}\in W_1$ and
\be
f^{(0)}(x-1)=if^{(0)}(x), \quad
f^{(1)}(x-1)=-if^{(1)}(x).
\label{16}
\ee
For $k=0,1,2,3$ set
\be
g^{(k)}(x)=\frac{1}{2}\sum_{l=0}^3\exm{\frac{l(1+4k)}{16}}f^{(0)}_{2,l/4}(x)=\sum_{l=0}^3\exm{\frac{l(1+4k)}{16}}
f^{(0)}\lll \frac {x-l}{4}\rrr.
\label{18}
\ee
We note that $g^{(k)}\in W_3$. Moreover, using~(\ref{16}) one can easily check that
\be
g^{(k)}(x-1)=\ex{ \frac {3+4k}{16} }g^{(k)}(x).
\label{17}
\ee
The vector-function $(f^{(0)}, g^{(0)}, g^{(1)}, g^{(2)}, g^{(3)})^T$ generates an ONWB because from (\ref{18})
it follows that it is unitary equivalent to
$(f^{(0)}, f_{2,0}^{(1)}, f_{2,1/4}^{(1)}, f_{2,1/2}^{(1)}, f_{2,3/4}^{(1)})^T$
and by Lemma~\ref{l2} the latter generates the same ONWB as $(f^{(0)}, f^{(1)})^T$.

Let $h^{(0)}=af^{(0)}+bg^{(2)}+cg^{(3)}$ where $a,b,c\in\bC$. Set
$$
h^{(1)}(x)=h^{(0)}(x-1).
$$
Then $h^{(1)}=\zeta^4 af^{(0)}+\zeta^{11}bg^{(2)}+\zeta^{15}cg^{(3)}$ by (\ref{16}) and (\ref{17})
where $\zeta=\ex{/16}$. So
$$
\langle h^{(0)}, h^{(1)}\rangle=\zeta^{12}|a|^2+\zeta^5|b|^2+\zeta|c|^2
$$
Take $a,b,c$ so that $|a|^2+|b|^2+|c|^2=1$ and the right-hand side of the last equality equal 0
(e.g., $a=\lambda$, $b={\lambda}{\sqrt{\cos \pi/8}}$,
$c={\lambda}{\sqrt{\sin\pi/8}}$, where $\lambda $ is a positive real). Then  there exists
a unitary matrix of the form
$$
U=\left( \begin{array}{rrr} a&b&c \\ \zeta^4a& \zeta^{11}b& \zeta^{15}c\\
\alpha&\beta&\gamma
\end{array} \right).
$$
Set $h^{(3)}=\alpha f^{(0)}+\beta g^{(2)}+\gamma g^{(3)}$. Then $(h^{(0)},h^{(1)}, h^{(2)})^T=U(f^{(0)}, g^{(2)}, g^{(3)})^T$.
So the vector function $(g^{(0)},g^{(1)}, h^{(0)}, h^{(1)}, h^{(2)})^T$ generates an ONWB. Finally, set
$$
h(x)=\frac{\sqrt2}2 h^{(0)}(2x).
$$
Then $h^{(0)}=h_{1,0}$, $h^{(1)}=h_{1,1/2}$, and consequently  by Lemma~\ref{l2} the vector function
$$
\Psi=(g^{(0)}, g^{(1)}, h^{(2)}, h)^T
$$
also generates an ONWB. We observe that $g^{(0)}, g^{(1)}, h^{(2)}\in V_4$ and $h\in V_3$.
Moreover, since $(\beta, \gamma)\ne(0, 0)$, the $W_3$-parts of the functions $g^{(0)}, g^{(1)}, h^{(2)}$
with respect to decomposition (\ref{3}) span a linear space of dimension 3.


Suppose that $\Psi=\Psi_0$ is reduced to a standard Haar vector-function $(\psi)=\Psi_N$ in $N$ steps.
There are two kind of steps: unitary ones for which $\rk(\Psi_j)=\rk(\Psi_{j-1})$,
and decreasing ones for which $\rk(\Psi_j)=\rk(\Psi_{j-1})-1$. Without loss of generality it can be assumed
that they alternate, the first step is unitary and the last is decreasing. Since $\rk(\Psi_0)=4$
and $\rk(\Psi_N)=1$, we have $N=6$, $\rk(\Psi_1)=4$ and  $\rk(\Psi_2)=3$. Moreover, $\Psi_1$ is in $V_4$,
$\Psi_2$ is in $V_3$ and the $W_3$-parts of components of $\Psi_1$ span a linear space of dimension
at most 2,
by the definition of decreasing step. So  $\Psi_0$ as unitary equivalent to $\Psi_1$, also has this property.
However, the dimension of the corresponding space for $\Psi$ equals 3 by above, a contradiction.
$\Diamond$

\bibliographystyle{amsplain}

\end{document}